# Synthèse d'observateur robuste pour la détection de défauts des systèmes linéaires discrets à commutation

Djamel. E. C. Belkhiat, Nadhir Messai, Noureddine Manamanni

Université de Reims Champagne Ardenne, CReSTIC
UFR SEN, Moulin de la Housse BP 1039, 51687, Reims cedex2, France.

prénom.nom@univ-reims.fr

*Résumé* —L'efficacité de l'opération de détection de défauts dans un système est tributaire des qualités du générateur de résidus utilisé. Dans ce cadre, ce papier traite de la synthèse d'observateurs hybrides robustes pour les systèmes linéaires discrets à commutations (SLDC). L'observateur est synthétisé afin de générer des résidus permettant de détecter les défauts affectant un SLDC. Ce dernier est supposé soumis à des effets indésirables dus à des entrées inconnues (bruit de mesure, perturbations externes, …) et/ ou des erreurs de modélisation. Pour ce faire, le problème de synthèse d'observateur est résolu à l'aide des outils de la commande robuste. Des conditions de convergence sous forme d'Inégalité Matricielles Linéaires (LMI) seront alors établies et une procédure itérative permettant de résoudre numériquement l'ensemble des LMI sera également proposée. Des résultats de simulation illustreront les performances de l'approche proposée.

*Mots clés* — Système Hybride, Système à commutations, Détection de défauts, $H_\infty$, Inégalités Matricielles Linéaires (LMI).

## I. INTRODUCTION

En vertu de l'essor fulgurant de l'étude des Systèmes Dynamiques Hybrides (SDH) au cours des deux dernières décennies, différents travaux portant sur la conception et la synthèse de loi de commande ont vu le jour [1] [2]. Toutefois, ces lois de commande deviennent souvent inefficaces dès qu'un défaut apparaît. Il est donc vital de détecter et d'isoler les défauts le plus rapidement possible afin d'éviter des comportements qui pourraient endommager le système au cours de son fonctionnement. Dans ce contexte, une attention particulière a été portée ces dernières années aux techniques de détection et d'isolation de défauts appliquées aux SDH [3] [4]. Ainsi, plusieurs approches de détection de défauts basées sur l'estimation d'état et/ ou l'identification paramétrique ont été adaptées aux SDH [5]. Cependant, peu de travaux de la littérature traitent le cas d'un SDH caractérisé par des erreurs de modélisation et des perturbations externes.

En effet, comme il n'est pas toujours évident de déterminer le modèle exact d'un système, les modèles utilisés sont souvent sujets à des erreurs de modélisation. Par ailleurs, les perturbations externes sont omniprésentes et inévitables particulièrement quand le système étudié renferme des composants électroniques sensibles. Deux qualités essentielles sont alors exigées dans un dispositif de détection et d'isolation de défauts. La première est sa robustesse vis-à-vis des perturbations externes et des erreurs de modélisation. Quant à la deuxième, elle concerne sa sensibilité aux défauts que l'on cherche à détecter [6].

Dans ce contexte, une stratégie robuste de détection de défaut pour une classe de SDH discret et incertain est présentée dans [7]. Cette approche est basée sur l'utilisation d'un observateur hybride qui génère des résidus qui sont à la fois robustes aux perturbations externes et sensibles aux défauts. Cette approche a été, ensuite, étendue dans [8] à une classe SDH non-linéaire discret et incertain et cela, en exploitant les résultats de [9] et de [10].

L'objectif de ce papier est de proposer une approche permettant de garantir un bon compromis entre la robustesse et la sensibilité d'un dispositif de détection et d'isolation de défauts pour une classe de systèmes linéaires discrets à commutations (SLDC). Ainsi, ce travail propose une approche de synthèse d'observateur pour cette classe de SLDC caractérisée par la présence de sauts d'états aux instants de commutations. Cet observateur aura pour rôle d'améliorer la tâche de détection de défauts en l'absence totale de connaissances sur le mode actif du SLDC. Pour ce faire, le problème est formulé en exploitant une approche d'optimisation $H_\infty$ [11], [12].

Cet article est organisé de la façon suivante : la section II présente la classe de SDH considérée dans ce travail. Notre démarche de synthèse d'observateur est donnée dans la section III. Finalement, un exemple de simulation illustrant l'efficacité de l'approche proposée est présenté dans la section IV.

## II. Système linéaire discret à commutations

La classe de système linéaire discret à commutations considérée dans ce travail est donnée par :

$$x(k+1) = A_{q(k)}x(k) + Bu(k) + B_{d_{q(k)}}d(k) + B_{f_{q(k)}}f(k) \quad (1)$$

$$y(k) = Cx(k) + D_{d_{q(k)}}d(k) + D_{f_{q(k)}}f(k) \quad (2)$$

où $x \in \mathfrak{R}^n$, $u \in \mathfrak{R}^m$, $y \in \mathfrak{R}^p$, $d \in \mathfrak{R}^m$ et $f \in \mathfrak{R}^l$ sont respectivement le vecteur d'état, le vecteur des entrées, le vecteur de sortie, le vecteur des entrées inconnues (perturbations, incertitudes non structurées, bruits) et le vecteur des défauts à détecter et à localiser. $A_{q(k)}, B, C, B_{f_{q(k)}}, B_{d_{q(k)}}, D_{f_{q(k)}}, D_{d_{q(k)}}$ sont des matrices connues de dimensions appropriées, $q(k) \in Q = 1, 2, ... N$ est l'état discret (indice du mode à l'instant $k$).

Dans ce travail, nous considérons que le passage d'un sous système à un autre ne se produit qu'à des endroits bien précis de l'espace d'état, noté $S_{i,j}$, qui sont décrits par des hyper plans linéaires :

$$S_{i,j} = \{x \in \mathfrak{R}^n \mid s_{i,j}x = 0\} \qquad (i,j) \in I_s \qquad (3)$$

Où : $I_s$ est un ensemble de uplets regroupant les différentes transitions possibles, entre deux modes, qui peuvent avoir lieu dans un système à commutations. On suppose également que la trajectoire discrète est connue à l'avance. Par contre, les instants de commutations ne le sont pas.

Aussi, dans la suite du papier, nous supposerons que le SLDC possède un nombre fini de changements de mode en un temps fini, et que tous les couples $(A_{q(k)}, C)$ sont observables. Enfin, sans perte de généralités, on suppose que les deux vecteurs $d, f$ sont bornés en norme $L_2$.

## III. SYNTHESE DE L'OBSERVATEUR ROBUSTE

La structure de l'observateur proposé est donnée par :

$$\hat{x}(k+1) = A_{\hat{q}(k)}\hat{x}(k) + Bu(k) + K_{\hat{q}(k)}(y(k) - \hat{y}(k)) \qquad (4)$$

$$\hat{y}(k) = C\hat{x}(k) \qquad (5)$$

$$r(k) = V(y(k) - \hat{y}(k)) \qquad (6)$$

où $\hat{x} \in \mathfrak{R}^n$ et $\hat{y} \in \mathfrak{R}^p$ sont respectivement les vecteurs d'état et de sortie estimés, $r$ est le vecteur des résidus, $K_{\hat{q}}$ et $V$ sont respectivement la matrice des gains et de pondération à déterminer.

Si on désigne par $q$ le mode du SLC et par $\hat{q}$ le mode de l'observateur, les dynamiques de l'erreur d'estimation $e = x - \hat{x}$ et du résidu $r$ sont données par :

$$e(k+1) = \bar{A}_{\hat{q}}e(k) + \Delta A_{q,\hat{q}}x(k) + \bar{B}_{d_q}d(k) + \bar{B}_{f_q}f(k) \qquad (7)$$

$$\begin{aligned} r(k) &= V(y(k) - \hat{y}(k)) \\ &= V(Ce(k) + D_{d_q}d(k) + D_{f_q}f(k)) \end{aligned} \qquad (8)$$

avec : $\bar{A}_{\hat{q}} = A_{\hat{q}} - K_{\hat{q}}C$, $\bar{B}_{d_q} = B_{d_q} - K_{\hat{q}}D_{d_q}$, $\bar{B}_{f_q} = B_{f_q} - K_{\hat{q}}D_{f_q}$, $\Delta A_{q,\hat{q}} = A_q - A_{\hat{q}}$.

On peut remarquer que la dynamique du résidu ne dépend pas seulement de $f$ et $d$ mais également de $x$. Ainsi, le problème de la robustesse de la détection de défauts revient à déterminer $K_{\hat{q}}$ et $V$ permettant d'avoir un résidu pertinent. Ce résidu devrait à la fois être sensible aux défauts $f$ et robuste aux entrées inconnues $d$. Par ailleurs, la présence des sauts d'états aux instants de commutation nécessite une actualisation de l'état aux instants de commutations (voir [13] pour plus détails).

La dynamique du résidu peut être alors reformulée, selon le principe de superposition des systèmes linéaires, de la façon suivante :

$$r(k) = r_x(k) + r_d(k) + r_f(k) \qquad (9)$$

Cette reformulation montre que la dynamique du résidu (9) est la somme de trois composantes que l'on peut traiter séparément. Ainsi, dans la suite, nous proposerons des conditions assurant à la fois la bornitude de la composante $r_x$ et la stabilité asymptotique des deux composantes $r_d$ et $r_f$.

Aussi, afin d'assurer la robustesse du résidu, nous proposons de minimiser le transfert des entrées inconnues $d$ sur le résidu $r$ selon l'inégalité suivante :

$$\|r_d\|_2^2 - \gamma_{\hat{q}}^2\|d\|_2^2 < 0 \qquad (10)$$

et de maximiser le transfert des défauts $f$ sur le résidu $d$ selon l'inégalité suivante :

$$\|r_f\|_2^2 - \beta_{\hat{q}}^2\|f\|_2^2 > 0 \qquad (11)$$

Avec les scalaires $\gamma_{\hat{q}}$ et $\beta_{\hat{q}}$ définis strictement positifs.

### 3.1 Synthèse de l'observateur

Avant de présenter le résultat principal de ce travail, nous présentons d'abord quelques outils nécessaires pour appréhender la suite.

**Définition 1** [14] :

La suite $(x(0), x(1), x(2), x(3)...)$ est dite bornée par $x_{\max}$ si $\forall k \geq 0 \quad \|x_k\| \leq x_{\max}$.

La suite $(x(0), x(1), x(2), x(3)...)$ est dite éventuellement bornée par $x_{\max}$ si $\forall \delta \geq 0 \; \exists K_0 > 0 \; \forall k \geq K_0 \; \|x(k)\| = x_{\max} + \delta$.

i.e. $\lim_{k \to \infty} \|x(k)\| \leq x_{\max}$.

**Lemme 1** (Lemme de Finsler's) [15]:

Etant donnés un vecteur $x \in \mathfrak{R}^n$ ainsi que les matrices $Q = Q^T \in \mathfrak{R}^{n \times n}$, $R \in \mathfrak{R}^{m \times n}$ et $S \in \mathfrak{R}^{m \times n}$ tels que $rang(R) < n$ et $rang(S) < n$ ; les expressions suivantes sont équivalentes:

- $x^TQx < 0 \quad \forall x \in \{x \in \mathfrak{R}^* \mid x \neq 0, Rx = 0, Sx = 0\}$.

- $R_\perp^T Q R_\perp < 0$ et $S_\perp^T Q S_\perp < 0$

- $\exists \sigma \in \mathfrak{R}$ tel que $Q - \sigma R^T R$ et $Q - \sigma S^T S < 0$.

- $\exists X \in \mathfrak{R}^{n \times m}$ tel que $Q + S^T XR + R^T X^T S < 0$.

où : $R_\perp$ est le complément orthogonal de $R$.

D'après les dynamiques (7) et (8), les relations suivantes sont déduites:

$$e_x(k+1) = \bar{A}_{\hat{q}}e_x(k) + \Delta A_{q,\hat{q}}x(k) \text{ et } r_x(k) = VCe_x(k).$$

A partir de là, nous pouvons énoncer le premier résultat de ce papier (Théorème 1). Ce dernier détermine les conditions garantissant la bornitude de la composante $e_x$ en présence de sauts d'états aux moments de commutations.

**Théorème 1.** Considérons le SLDC (1)-(2) avec l'observateur (4)-(6) et supposons que pour $K_0 > 0$ nous avons $\max_{k > K_0} \|x(k)\| \leq x_{\max}$.

Supposons également qu'il existe des matrices $P_{\hat{q}} = P_{\hat{q}}^T \geq 0$, $K_{\hat{q}}$, et des scalaires $\sigma_{\hat{q}}, \varepsilon > 0, \alpha > 0, \xi > 0, v_{\hat{q},q} \geq 0$ tels que les conditions I, II données ci-après sont vérifiées.

I. l'ensemble des inégalités matricielles suivant (12)-(14) est satisfait:

1. $\alpha I \leq P_{\hat{q}} \leq \xi I$, $\hat{q} \in I_N$ (12)

2. $\begin{bmatrix} \rho - \frac{1}{4}\sigma_{\hat{q}}C^TC + v_{\hat{q},q}I & \bar{A}_{\hat{q}}^T P_{\hat{q}}\Delta A_{q,\hat{q}} \\ (*) & \Delta A_{q,\hat{q}}^T P_{\hat{q}}\Delta A_{q,\hat{q}} - \varepsilon^2 v_{\hat{q},q}I \end{bmatrix} \leq 0$ (13)

3. $P_q = P_{\hat{q}} + d_{\hat{q},q}C + C^T d_{\hat{q},q}$, $(\hat{q},q) \in I_s$ (14)

avec: $\rho = A_{\hat{q}}^T P_{\hat{q}} A_{\hat{q}} - A_{\hat{q}}^T P_{\hat{q}} K_{\hat{q}} C - C^T K_{\hat{q}}^T P_{\hat{q}} A_{\hat{q}} - P_{\hat{q}}$.

II. $\forall \hat{x} \in S_{i,j}$ l'état de l'observateur hybride est actualisé selon:

$$\hat{x}^+ = \left(I - R_{\hat{q}}^{-1}\left(CR_{\hat{q}}^{-1}\right)^\dagger C\right)\hat{x} + R_{\hat{q}}^{-1}\left(CR_{\hat{q}}^{-1}\right)^\dagger y \quad (15)$$

Alors, l'erreur d'estimation $e_x$ est éventuellement bornée (au sens de la définition 1) par $e_{x\max} \leq \sqrt{\frac{\xi}{\alpha}} \varepsilon x_{\max}$.

*Preuve :* En introduisant la fonction de Lyapunov multiple.

$$V_{\hat{q}}\left(e_x(k)\right) = e_x^T(k) P_{\hat{q}} e_x(k) \qquad \hat{q} \in I_N \quad (16)$$

où : $P_{\hat{q}} \in \Re^{n \times n}$ sont des matrices symétriques définies positives. Alors, $\Delta V_{\hat{q}}\left(e_x(k)\right) = V_{\hat{q}}\left(e_x(k+1)\right) - V_{\hat{q}}\left(e_x(k)\right)$ est donnée par:

$$\begin{aligned}\Delta V_{\hat{q}}\left(e_x(k)\right) &= e_x^T(k) \bar{A}_{\hat{q}}^T P_{\hat{q}} \bar{A}_{\hat{q}} e_x(k) \\ &+ e_x^T(k) \bar{A}_{\hat{q}}^T P_{\hat{q}} \Delta A_{q,\hat{q}} x(k) \\ &+ x^T(k) \Delta A_{q,\hat{q}}^T P_{\hat{q}} \Delta A_{q,\hat{q}} x(k) \\ &+ x^T(k) \Delta A_{q,\hat{q}}^T P_{\hat{q}} \bar{A}_{\hat{q}} e_x(k) - e_x^T(k) P_{\hat{q}} e_x(k)\end{aligned} \quad (17)$$

L'équation (17) peut être écrite sous la forme d'inégalité matricielle garantissant la stabilité au sens de Lyaponov:

$$\begin{bmatrix} \rho + C^T K_{\hat{q}}^T P_{\hat{q}} K_{\hat{q}} C & \bar{A}_{\hat{q}}^T P_{\hat{q}} \Delta A_{q,\hat{q}} \\ (*) & \Delta A_{q,\hat{q}}^T P_{\hat{q}} \Delta A_{q,\hat{q}} \end{bmatrix} \leq 0 \quad (18)$$

où $\rho = A_{\hat{q}}^T P_{\hat{q}} A_{\hat{q}} - A_{\hat{q}}^T P_{\hat{q}} K_{\hat{q}} C - C^T K_{\hat{q}}^T P_{\hat{q}} A_{\hat{q}} - P_{\hat{q}}$.

Notons que l'inégalité (18) n'est pas encore sous forme LMI et elle ne peut être définie négative si $\Delta A_{\hat{q},q} \neq 0$. Ainsi, en exploitant la S-procédure pour la contrainte quadratique $\|e_x\|^2 \geq \varepsilon^2 \|x\|^2$, l'inégalité matricielle (19) peut être vérifiée facilement :

$$\begin{bmatrix} \rho + C^T K_{\hat{q}}^T P_{\hat{q}} K_{\hat{q}} + v_{\hat{q},q} I & \bar{A}_{\hat{q}}^T P_{\hat{q}} \Delta A_{q,\hat{q}} \\ (*) & \Delta A_{q,\hat{q}}^T P_{\hat{q}} \Delta A_{q,\hat{q}} - \varepsilon^2 v_{\hat{q},q} I \end{bmatrix} \leq 0 \quad (19)$$

L'inégalité (19) est ensuite développée comme suit :

$$\underbrace{\begin{bmatrix} \rho + v_{\hat{q},q} I & \bar{A}_{\hat{q}}^T P_{\hat{q}} \Delta A_{q,\hat{q}} \\ (*) & \Delta A_{q,\hat{q}}^T P_{\hat{q}} \Delta A_{q,\hat{q}} - \varepsilon^2 v_{\hat{q},q} I \end{bmatrix}}_{Q} \\ + \underbrace{\begin{bmatrix} \frac{1}{2} C^T \\ 0 \end{bmatrix}}_{R^T} M_{\hat{q}} \underbrace{\begin{bmatrix} C & 0 \end{bmatrix}}_{S} + \underbrace{\begin{bmatrix} C^T \\ 0 \end{bmatrix}}_{S^T} M_{\hat{q}} \underbrace{\begin{bmatrix} \frac{1}{2} C & 0 \end{bmatrix}}_{R} \leq 0 \quad (20)$$

avec $M_{\hat{q}} = K_{\hat{q}}^T P_{\hat{q}} K_{\hat{q}}$.

Le lemme 1, montre que l'inégalité matricielle (20) est équivalente aux inégalités matricielles (21) et (22):

$$\begin{bmatrix} \rho + v_{\hat{q},q} I & \bar{A}_{\hat{q}}^T P_{\hat{q}} \Delta A_{q,\hat{q}} \\ (*) & \Delta A_{q,\hat{q}}^T P_{\hat{q}} \Delta A_{q,\hat{q}} - \varepsilon^2 v_{\hat{q},q} I \end{bmatrix} - \sigma_{\hat{q}} \begin{bmatrix} C^T C & 0 \\ 0 & 0 \end{bmatrix} \leq 0 \quad (21)$$

$$\begin{bmatrix} \rho + v_{\hat{q},q} I & \bar{A}_{\hat{q}}^T P_{\hat{q}} \Delta A_{q,\hat{q}} \\ (*) & \Delta A_{q,\hat{q}}^T P_{\hat{q}} \Delta A_{q,\hat{q}} - \varepsilon^2 v_{\hat{q},q} I \end{bmatrix} - \frac{\sigma_{\hat{q}}}{4} \begin{bmatrix} C^T C & 0 \\ 0 & 0 \end{bmatrix} \leq 0 \quad (22)$$

On peut remarquer que :

$$\begin{bmatrix} \rho + v_{\hat{q},q} I & \bar{A}_{\hat{q}}^T P_{\hat{q}} \Delta A_{q,\hat{q}} \\ (*) & \Delta A_{q,\hat{q}}^T P_{\hat{q}} \Delta A_{q,\hat{q}} - \varepsilon^2 v_{\hat{q},q} I \end{bmatrix} \leq \frac{\sigma_{\hat{q}}}{4} \begin{bmatrix} C^T C & 0 \\ 0 & 0 \end{bmatrix} \\ \leq \sigma_{\hat{q}} \begin{bmatrix} C^T C & 0 \\ 0 & 0 \end{bmatrix} \quad (23)$$

Ceci nous permet de garder seulement l'inégalité (22), qui vérifie implicitement l'inégalité (21).

Concernant les démonstrations pour aboutir à (14) et (15), le lecteur peut se référer aux travaux dans [13] et les références citées dans son article. Ces deux conditions sont introduites dont l'objectif est d'assurer que l'énergie globale de l'erreur d'estimation $e_x(t)$ est décroissante aux instants de commutations tel que cette inégalité est vérifiée:

$$\left(x - \hat{x}^+\right)^T P_q \left(x - \hat{x}^+\right) \leq \left(x - \hat{x}\right)^T P_{\hat{q}} \left(x - \hat{x}\right)$$

Il s'agit maintenant de prouver que l'erreur d'estimation est éventuellement bornée (au sens de la définition 1) par une constante.

Considérons la contrainte quadratique suivante :

$$\|e_x\|^2 \geq \varepsilon^2 \|x\|^2 \quad (24)$$

Pour quelques valeurs de $\varepsilon > 0$, nous supposons que $\Delta V_{\hat{q}}\left(e_x(k)\right) < 0$ quand (24) est vérifiée.

Pour une valeur arbitraire de $\delta > 0$, on note par :

$$V_{\max}^\delta = \max_{\|e_x\| \leq \varepsilon x_{\max} + \delta} V_{\hat{q}}\left(e_x(k)\right) \quad (25)$$

On désigne par $S_\delta$ l'ensemble suivant:

$$S_\delta = \left\{e_x \mid V_{\hat{q}}\left(e_x(k)\right) < V_{\max}^\delta\right\} \quad (26)$$

Puisque $\Delta V_{\hat{q}}\left(e_x(k)\right) < 0$ pour $e_x \notin S_\delta$, il s'ensuit que $S_\delta$ est positivement invariant c'est à dire pour $K_0 > 0$.

$$\exists K_0 > 0 \quad V_{\hat{q}}\left(e_x(K_0)\right) < V_{\max}^\delta$$

Selon (12), il résulte que :

$$V_{\max}^\delta < \xi \left[\varepsilon x_{\max} + \delta\right]^2 \quad (27)$$

Par conséquent,

$$\forall_{\delta > 0} \exists_{K_0 > 0} \forall_{k > K_0} \quad \|e_x(k)\| \leq \sqrt{\frac{\xi}{\alpha}}\left[\varepsilon x_{\max} + \delta\right] \quad (28)$$

c'est à dire

$$e_{x\max} = \lim_{k \to \infty} \|e_x(k)\| \leq \sqrt{\frac{\xi}{\alpha}} \varepsilon x_{\max} \qquad \blacksquare$$

**Remarque 1.** Afin de résoudre les inégalités matricielles du théorème 1 avec des outils mathématiques classiques et aussi pour éviter que les gains d'observateur $K_{\hat{q}}$ soient trop grands, nous introduisons une nouvelle variable $H_{\hat{q}} \in \Re^{n \times p}$ et la condition suivante :

$$\begin{bmatrix} \lambda_{\hat{q}}^2 I_{p \times p} & * \\ H_{\hat{q}} & I_{n \times n} \end{bmatrix} \geq 0 \qquad \hat{q} \in I_N \quad (29)$$

où $\lambda_{\hat{q}}$ est un paramètre de conception et $H_{\hat{q}} = P_{\hat{q}} K_{\hat{q}}$.

### 3.2 Robustesse de l'observateur

Dans cette section, les deux composantes du résidu $r_x$ et $r_u$ seront traitées séparément. Ainsi, nous formulons le problème de la synthèse de l'observateur robuste comme étant un problème d'optimisation, puis, nous le résolvons en exploitant les résultats de la commande robuste.

D'après les dynamiques (7) et (8), les relations suivantes sont déduites :

$$\begin{cases} e_d(k+1) = \overline{A}_{\hat{q}} e_d(k) + \overline{B}_{d_q} d(k) \\ r_d(k) = VC e_d(k) + VD_{d_q} d(k) \end{cases} \quad (30)$$

et

$$\begin{cases} e_f(k+1) = \overline{A}_{\hat{q}} e_f(k) + \overline{B}_{f_q} f(k) \\ r_f(k) = VC e_f(k) + VD_{f_q} f(k) \end{cases} \quad (31)$$

**Lemme 2.** Considérons le SLC (1)-(2) avec l'observateur (4)-(6). Le système (30) est asymptotiquement stable et satisfait la condition (10) s'il existe des matrices $P_{\hat{q}} > 0$, $H_{\hat{q}}$ et $V$ et le scalaire $\kappa_{\hat{q}}$ tels que la LMI suivante est vérifiée :

$$\begin{bmatrix} \rho - \dfrac{\kappa_{\hat{q}}}{4} C^T C & \theta - \dfrac{\kappa_{\hat{q}}}{4} C^T D_{d_q} & C^T V^T \\ (*) & \lambda - \dfrac{\kappa_{\hat{q}}}{4} D_{d_q}^T D_{d_q} & D_{d_q}^T V^T \\ (*) & (*) & -I \end{bmatrix} \leq 0 \quad (32)$$

avec : $\rho = A_{\hat{q}}^T P_{\hat{q}} A_{\hat{q}} - A_{\hat{q}}^T H_{\hat{q}} C - C^T H_{\hat{q}}^T A_{\hat{q}} - P_{\hat{q}}$

$\theta = A_{\hat{q}}^T P_{\hat{q}} B_{d_q} - A_{\hat{q}}^T H_{\hat{q}} D_{d_q} - C^T H_{\hat{q}}^T B_{d_q}$,

$\lambda = B_{d_q}^T P_{\hat{q}} B_{d_q} - B_{d_q}^T H_{\hat{q}} D_{d_q} - D_{d_q}^T H_{\hat{q}}^T B_{d_q} - \gamma_{\hat{q}}^2 I$.

*Preuve :* En partant du système (30), la condition (10) peut être alors reformulée de la façon suivante:

$$\sum_{k=0}^{N} \left( r_d^T r_d \right) \leq \gamma_{\hat{q}}^2 \sum_{k=0}^{N} \left( d^T d \right)$$

En introduisant une fonction multiple de Lyapunov candidate $V(e_d) = e_d^T P_{\hat{q}} e_d \geq 0$, $P_{\hat{q}} > 0$, nous obtenons :

$$J = \sum_{k=0}^{N} \left( r_d^T r_d - \gamma_{\hat{q}}^2 d^T d + \Delta V(e_d(k)) \right) - V(e_d(N))$$

$$J = \left( \sum_{k=0}^{N} \begin{bmatrix} e_d^T & d^T \end{bmatrix} \underbrace{\left[ \begin{bmatrix} C^T V^T \\ D_{d_q}^T V^T \end{bmatrix} \begin{bmatrix} VC & VD_{d_q} \end{bmatrix} + E_d \right]}_{R_d} \begin{bmatrix} e_d \\ d \end{bmatrix} \right) - V(e(N))$$

où $E_d = \begin{bmatrix} \overline{A}_{\hat{q}}^T P_{\hat{q}} \overline{A}_{\hat{q}} - P_{\hat{q}} & \overline{A}_{\hat{q}}^T P_{\hat{q}} \overline{B}_{d_q} \\ (*) & \overline{B}_{d_q}^T P_{\hat{q}} \overline{B}_{d_q} - \gamma_{\hat{q}}^2 I \end{bmatrix}$. D'où, $J \leq 0$, si $R_d \leq 0$.

En utilisant le complément de Schur, la matrice $R_d$ devient :

$$R_d = \begin{bmatrix} \rho + C^T K_{\hat{q}}^T P_{\hat{q}} K_{\hat{q}} C & \theta_d + C^T K_{\hat{q}}^T P_{\hat{q}} K_{\hat{q}} D_{d_q} & C^T V^T \\ (*) & \lambda_d + D_{d_q}^T K_{\hat{q}}^T P_{\hat{q}} K_{\hat{q}} D_{d_q} & D_{d_q}^T V^T \\ (*) & (*) & -I \end{bmatrix} \leq 0 \quad (33)$$

Comme pour le théorème 1, on propose de développer l'inégalité matricielle (33) à l'aide du lemme 1. Le reste du développement sera omis dans la suite vu le nombre limité de page. ∎

**Lemme 3.** Considérons le SLC (1)-(2) avec l'observateur (4)-(6). Le système (31) est asymptotiquement stable et satisfait la condition (11) s'il existe des matrices $P_{\hat{q}} > 0$, $H_{\hat{q}}$ et $V$ et le scalaire $\varsigma_{\hat{q}}$ tels que la LMI suivante est vérifiée :

$$\begin{bmatrix} \rho - \dfrac{\varsigma_{\hat{q}}}{4} C^T C + 2\varphi_1(V, V_c^n) & -\theta_f + \dfrac{\varsigma_{\hat{q}}}{4} C^T D_{f_q} & C^T V^T \\ (*) & \lambda_f - \dfrac{\varsigma_{\hat{q}}}{4} D_{f_q}^T D_{f_q} + 2\varphi_2(V, V_f^n) & D_{f_q}^T V^T \\ (*) & (*) & -I \end{bmatrix} \leq 0 \quad (34)$$

Avec: $\lambda_f = B_{f_q}^T P_{\hat{q}} B_{f_q} - A_{\hat{q}}^T H_{\hat{q}} D_{f_q} - D_{f_q}^T H_{\hat{q}}^T B_{f_q} + \beta_{\hat{q}}^2 I$,

$\theta_f = A_{\hat{q}}^T P_{\hat{q}} B_{f_q} - A_{\hat{q}}^T H_{\hat{q}} D_{f_q} - C^T H_{\hat{q}}^T B_{f_q}$,

$\varphi_1(V, V_c^n) = (V_c^n)^T V_c^n - (V_c^n)^T VC - C^T V^T V_c^n$,

$\varphi_2(V, V_f^n) = (V_f^n)^T V_f^n - (V_f^n)^T V D_{f_q} - D_{f_q}^T V^T V_f^n$.

*Preuve :* de la même manière que précédemment, la condition (11) est reformulée de la façon suivante:

$$\sum_{k=0}^{N} \left( r_f^T r_f \right) \geq \beta_{\hat{q}}^2 \sum_{k=0}^{N} \left( f^T f \right)$$

En introduisant une fonction multiple de Lyapunov candidate $V(e_f) = e_f^T P_{\hat{q}} e_f \geq 0$, $P_{\hat{q}} > 0$, nous obtenons :

$$J = \left( \sum_{k=0}^{N} \begin{bmatrix} e_f^T & f^T \end{bmatrix} \underbrace{\left[ \begin{bmatrix} C^T V^T \\ D_{f_q}^T V^T \end{bmatrix} \begin{bmatrix} VC & VD_{f_q} \end{bmatrix} - E_f \right]}_{R_f} \begin{bmatrix} e_f \\ f \end{bmatrix} \right) + V(e(N))$$

avec $E_f = \begin{bmatrix} \overline{A}_{\hat{q}}^T P_{\hat{q}} \overline{A}_{\hat{q}} - P_{\hat{q}} & \overline{A}_{\hat{q}}^T P_{\hat{q}} \overline{B}_{f_q} \\ (*) & \beta_{\hat{q}}^T I + \overline{B}_{f_q}^T P_{\hat{q}} \overline{B}_{f_q} \end{bmatrix}$

ainsi, si

$$R_f = \begin{bmatrix} C^T V^T V C - \overline{A}_{\hat{q}}^T P_{\hat{q}} \overline{A}_{\hat{q}} + P_{\hat{q}} & C^T V^T V D_{f_q} - \overline{A}_{\hat{q}}^T P_{\hat{q}} \overline{B}_{f_q}^T \\ (*) & D_{f_q}^T V^T V D_{f_q} - \beta_{\hat{q}}^2 I - \overline{B}_{f_q}^T P_{\hat{q}} \overline{B}_{f_q} \end{bmatrix} \geq 0, \quad (35)$$

on peut donc obtenir $J \geq 0$.

L'inégalité (35) est malheureusement conservative, vu la présence des termes non-linéaires $V^T V$, $P_{\hat{q}} K_{\hat{q}}$ et $K_{\hat{q}}^T P_{\hat{q}} K_{\hat{q}}$, et ne peut être résolue à l'aide des outils mathématiques classiques. C'est pour ces raisons que l'on propose le développement de l'inégalité (35) à l'aide du Lemme 1 et le complément de Schur. ∎

Les résultats précédents nous permettent maintenant d'énoncer le théorème 2.

**Théorème 2.** Considérons le SLDC (1)-(2) avec l'observateur (4)-(6), le résidu (8) est borné et satisfait (10) et (11), si :

1. il existe des matrices $P_{\hat{q}} > 0$, $H_{\hat{q}}$ et $V$ et des scalaires $\sigma_{\hat{q}}$, $\kappa_{\hat{q}}$, $\varsigma_{\hat{q}}$ tels que les LMIs (12), (13), (14), (29), (32) et (34) soient vérifiées.
2. l'état estimé est actualisé selon (15) à l'instant de commutation.

*Preuve :* D'après les propriétés algébriques de la somme, chaque composante de la dynamique de résidus peut être traitée séparément.
1. Conformément au théorème 1 (LMIs (12), (13), (14) et (29)) la première composante du résidu (implicitement, de l'erreur d'estimation) est bornée.
2. conformément au lemme 2 (LMI (32)), la seconde composante du résidu (implicitement, de l'erreur d'estimation) est asymptotiquement stable et satisfait la condition (10).
3. conformément au lemme 3 (LMI (34)), la troisième composante du résidu (implicitement, de l'erreur d'estimation) est asymptotiquement stable et satisfait la condition (11).

Finalement, en combinant les trois points 1, 2 et 3, on peut facilement prouver le théorème 2. ∎

### 3.3 Optimisation de la Robustesse vs sensibilité

Les valeurs de $P_{\hat{q}} > 0$, $H_{\hat{q}}$ et $V$ satisfaisant le théorème 2 et minimisant le critère $Indice\_Perf = \inf_{\omega} \frac{\gamma_{\hat{q}}}{\beta_{\hat{q}}}$, peuvent être obtenus en utilisant la procédure itérative suivante.

1. Initialiser $\gamma_{\hat{q}}$ et $\beta_{\hat{q}}$.
2. résoudre les LMIs (12), (13), (14), (29) et (32) afin de trouver une solution faisable pour les matrices $P_{\hat{q}}$, $H_{\hat{q}}$ et $V$.
3. si les LMIs (12), (13), (14), (29) et (32) ne sont pas faisables, prendre d'autres valeurs pour $\gamma_{\hat{q}}$ et $\beta_{\hat{q}}$, et refaire l'étape jusqu'à ce qu'une solution soit trouver.
4. calculer $V^{n-1} = V$, $V_c^n = V^{n-1} C_{\hat{q}}$, $V_f^n = V^{n-1} F_{f_q}$ et substituer $P_{\hat{q}}$ par sa valeur trouvée à l'étape 2 dans les LMIs (12), (13), (14), (29), (32) et (34) afin de trouver une autre solution faisable pour $H_{\hat{q}}$ et $V$ (i.e., dans cette étape, le terme $P_{\hat{q}}$ n'est plus une matrice variable).
5. décrémenter $\gamma_{\hat{q}}$ et incrémenter $\beta_{\hat{q}}$ puis répéter l'étape 4 jusqu'à ce qu'aucune solution faisable ne soit obtenue.
6. Récupérer les dernières valeurs de $H_{\hat{q}}$ et $V$.

### 3.4 Détection de défauts

La logique de détection considérée dans ce travail est donnée par :

$\|r\|_{2,K} > J_{th}$, alarme, défaut est détecté. (36)

$\|r\|_{2,K} \leq J_{th}$, aucune alarme, aucun défaut n'est détecté (37)

où $\|r\|_{2,K}$ est déterminée :

$$\|r\|_{2,K} = \left[\sum_{k_1}^{k_2} r(k)^T r(k)\right]^{\frac{1}{2}} = \|r_x(k) + r_d(k) + r_f(k)\|_{2,K}, \quad K = k_2 - k_1$$

avec $J_{th} = J_{th,d} = \sup_{x \in L_2, d \in L_2} \|r_c(k)\|_{2,K}$, $r_c(k) = r(k)|_{f=0}$ et $r_f(k) = r(k)|_{x=0, d=0}$.

## IV. SIMULATION ET RESULTATS

On considère, dans cet exemple, le système à commutations à deux modes défini comme suit :

Mode 1, si $S_1 = 0$ : $A_1 = \begin{bmatrix} 1 & 2*10^{-4} \\ -2*10^{-4} & 1 \end{bmatrix}$, $B = \begin{bmatrix} 0 \\ 0 \end{bmatrix}$,

$C = \begin{bmatrix} 1.01 & -2.4 \end{bmatrix}$, $B_{d_1} = \begin{bmatrix} 1 & 0 \\ 1 & 0 \end{bmatrix}$, $B_{f_1} = \begin{bmatrix} 0 \\ 0 \end{bmatrix}$, $D_{d_1} = \begin{bmatrix} 1 & 0 \end{bmatrix}$,

$D_{f_1} = 0$.

Mode 2, si $S_2 = 0$ : $A_2 = \begin{bmatrix} 1 & 10^{-3} \\ 2*10^{-4} & 1 \end{bmatrix}$, $B = \begin{bmatrix} 0 \\ 0 \end{bmatrix}$,

$C = \begin{bmatrix} 1.01 & -2.4 \end{bmatrix}$, $B_{d2} = \begin{bmatrix} 1 & 0 \\ 0 & 1 \end{bmatrix}$, $B_{f_2} = \begin{bmatrix} 0 \\ 0 \end{bmatrix}$, $D_{d_2} = \begin{bmatrix} 1 & 1 \end{bmatrix}$,

$D_{f_2} = 4.3$.

où $S_i = a_i x_1 + b_i x_2, i \in \{1, 2\}$ sont les lois de commutations avec $a_1 = 0.25$, $b_1 = 1$, $a_2 = 1$ et $b_2 = 1$.

On note que le système à commutations est stable comme montre la figure1.

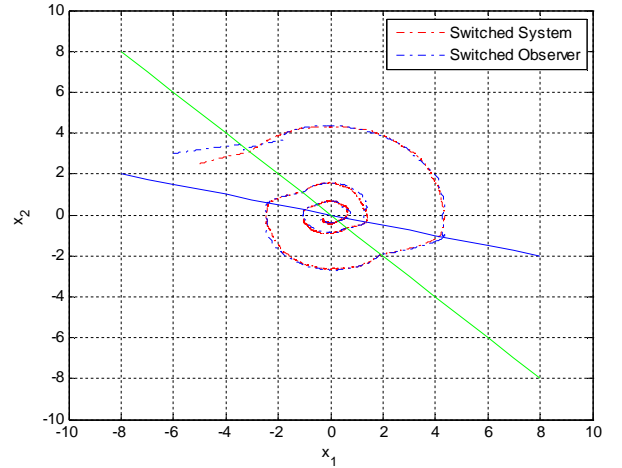

Fig.1. Plan de phase.

L'application de l'approche proposée pour les valeurs suivantes $\xi = 1.0193$, $\alpha = 0.0065$ et $\varepsilon = 2.3975$, permet d'aboutir a :

$\gamma_1 = \gamma_2 = 1.7398$, $\beta_1 = \beta_2 = 2.433$, $\frac{\gamma_{\hat{q}}}{\beta_{\hat{q}}} = 0.7150$ et $\lambda_1 = \lambda_2 = 20.044$.

$K_1 = \begin{bmatrix} 3.5616 \\ -0.8557 \end{bmatrix}$, $K_2 = \begin{bmatrix} -0.0434 \\ 0.0545 \end{bmatrix}$ et $V = 0.0174$.

Afin de montrer l'efficacité de l'approche proposée, les entrées inconnues $d$ sont considérées comme un bruit blanc gaussien centré d'une amplitude de 0.025 (Temps d'échantillonnage 0.001s). Le défaut $f$ est simulé comme une impulsion, d'une amplitude de 0.02, qui se produit en mode 2, entre les instants 13.5s et 13.8s. Le résidu généré est illustré par la figure 2, qui montre que le résidu converge vers zéro à l'instant 4s, ensuite, il s'éloigne notablement de zéro à la suite de l'apparition de défauts.

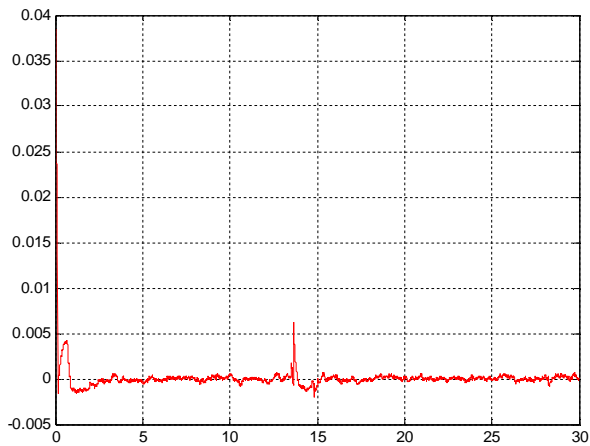

Fig.2. Résidu généré par RFDO.

De plus, la figure 3 illustre l'évolution de la fonction $\|r\|_{2,T}$ pour $T = 0.1s$. On peut remarquer que le défaut est détecté à l'instant 13.62s qui correspond à un retard de détection 0.12s.

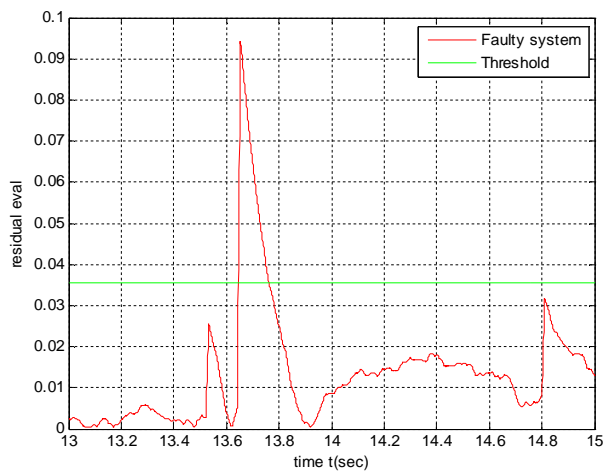

Fig.3. détection de défauts.

## V. Conclusion

Dans cet article, un observateur hybride robuste est synthétisé pour une classe de système linéaire discret à commutations. Cet observateur a pour vocation essentielle de détecter les défauts d'un SLDC soumis aux effets nuisibles de perturbations et/ou d'erreurs de modélisation. L'approche proposée est basée sur l'utilisation de résultats développés en matière de commande robuste. Finalement un exemple numérique illustre l'efficacité de l'approche développée.